\documentclass[12pt]{amsart}
\usepackage{times}
\usepackage{amsfonts}
\usepackage{amsthm}
\usepackage{amsmath}
\advance \topmargin by -\headheight
\advance \topmargin by -\headsep
\evensidemargin \oddsidemargin
\DeclareMathOperator{\tr}{tr}

\newtheorem{theorem}{Theorem}[section]
\newtheorem*{result}{Theorem}
\newtheorem{lemma}[theorem]{Lemma}

\newtheorem{remark}[theorem]{Remark}

\newtheorem{notation}[theorem]{Notation}

\begin{document}

\title[Free dimension and property $T$]{All generating sets of 
all property T von Neumann algebras\\ have 
free entropy dimension $\leq 1$}

\begin{abstract} Suppose that $N$ is a diffuse, property T von Neumann algebra 
and $X$ is an arbitrary finite generating set of selfadjoint elements for 
$N$.  By using rigidity/deformation arguments applied to 
representations of $N$ in ultraproducts of full matrix algebras, we deduce 
that the microstate spaces of $X$ are asymptotically discrete up to 
unitary conjugacy.  We use this description to show that the free entropy 
dimension of $X$, $\delta_0(X)$, is less than or equal to $1$.  It follows 
that when $N$ embeds into the ultraproduct of the hyperfinite 
$\mathrm{II}_1$-factor, then $\delta_0(X) =1$ and otherwise, $\delta_0(X) 
= -\infty$. This generalizes the earlier results of Voiculescu,
and Ge, Shen pertaining to $SL_n(\mathbb Z)$ as well as the results
of Connes, Shlyakhtenko pertaining to group generators of arbitrary
property $T$ algebras.

\end{abstract} \author{Kenley Jung \and Dimitri Shlyakhtenko}

\address{Department of Mathematics, University of California,
Los Angeles, CA 90095-1555,USA}

\email{kjung@math.ucla.edu, shlyakht@math.ucla.edu}
\thanks{Research supported in part by the NSF}
\maketitle

\section*{Introduction} In \cite{dvv:entropy2} and \cite{dvv:entropy3}, Voiculescu introduced the notion 
of free entropy dimension. For $X$ a finite set of self-adjoint elements 
of a tracial von Neumann algebra, $\delta_0(X)$ is a kind of asymptotic 
Minkowski dimension of the set of matricial microstates for $X$.  These 
notions led to the solution of several old operator algebra problems 
(see \cite{dvv:entropysurvey} for an overview).  Closely tied to this is the 
invariance question for $\delta_0$ which asks the following.  If $X$ and 
$Y$ are two finite sets of selfadjoint elements generating the same tracial von Neumann 
algebra, then is it true that 
$\delta_0(X) = \delta_0(Y)$?

For certain $X$ one can compute $\delta_0(X)$ and answer the invariance 
question in the affirmative.  
Suppose that $N=W^*(X)$ is diffuse and 
embeds into the ultraproduct of the hyperfinite $\mathrm{II}_1$-factor.  
Then $\delta_0(X) =1$ when $N$ has property $\Gamma$, or has a Cartan 
subalgebra, or is nonprime, or can be decomposed as an amalgamated free 
product of these algebras over a common diffuse subalgebra (see 
\cite{ge-shen:commutingGenerators,jung-freexentropy,jung:onebounded,dvv:entropy3}).

Another class of algebras to investigate in regard to possible values of 
$\delta_0(X)$ and the invariance question are those with Kazhdan's 
property T (\cite{connes:propertyT,kazhdan:T,popa:correspondences}).  These 
first appeared in the von Neumann algebra context in Connes' seminal work 
\cite{connes:fundgroup}. In recent years, Popa introduced the technique of 
playing the rigidity properties of such algebras against 
deformation results; this has led to a number of significant advances in 
the theory of von Neumann algebras.
 (\cite{popa:correspondences}, \cite{popa:betti}, 
\cite{ioana-peterson-popa:noOut}). 

Voiculescu made the first computations of $\delta_0$ for property T 
factors by showing that if $x_1, \ldots, x_n$ are diffuse, selfadjoint 
elements in a tracial von Neumann algebra such that for each $1 \leq i 
\leq n-1$, $x_i x_{i+1} = x_{i+1} x_i$, then $\delta_0(x_1,\ldots x_n) 
\leq 1$.  For $n \geq 3$, there exists a finite set of generators $X_n$ 
for the group algebra $\mathbb C SL_n(\mathbb Z)$ with this property
(this was first used in the context of measurable equivalence 
relations by Gaboriau \cite{gaboriau:cost} to prove that their 
cost is at most $1$). 
Hence $L(SL_n(\mathbb Z))$ has a set of generators $X$ for which 
$\delta_0(X) \leq 1$.  This was generalized  in
\cite{ge-shen:commutingGenerators} (see also 
\cite{ge:icmsurvey} and references therein) where Ge and Shen weakened the 
conditions on the generators $x_i$ and in particular obtained the stronger 
statement that $\delta_0(Y)\leq 1$ for any other set $Y$ of self-adjoint 
generators of the  von Neumann algebra.
However, all of these results rely on the special 
algebraic properties of certain generators (e.g. in $SL_n(\mathbb Z)$) and thus 
do not apply to the more general property $T$ groups or von Neumann 
algebras.

In \cite{connes-shlyakht:l2betti} a notion of $L^2$ 
cohomology for von Neumann algebras was introduced, and the values of the resulting
$L^2$ Betti numbers were connected with free probability 
and  the value of $\delta_0$.  Indeed, using cohomological ideas, it was 
proved in \cite{connes-shlyakht:l2betti} that if $X\subset 
\mathbb C\Gamma$ is an arbitrary set of generators, then $$ 
\delta_0(X)\leq \beta_1^{(2)}(\Gamma)-\beta_0^{(2)}(\Gamma)+1.$$ Here 
$\beta_j^{(2)}(\Gamma)$ are the Atiyah-Cheeger-Gromov $\ell^2$-Betti 
numbers of $\Gamma$ (see e.g. \cite{luck:book}).  This inequality is quite 
complicated to prove; indeed, one first proves the same inequality with 
$\delta_0$ replaced by its ``non-microstates'' analog $\delta^*$, and then 
uses a highly nontrivial result of Biane, Capitaine, Guionnet 
\cite{guionnet-biane-capitaine:largedeviations} that implies $\delta_0\leq 
\delta^*$.

In the case that $\Gamma$ has property T, the first $\ell^2$ Betti number 
vanishes (\cite{cheeger-gromov:l2}). So for $\Gamma$ an infinite group, one 
has $\delta_0(X)\leq 1$ for any finite generating set $X\subset \mathbb 
C\Gamma$.  However, even in this case, an ``elementary'' proof of this 
bound was not available and, moreover, it was not known whether 
$\delta_0(X) \leq 1$ for any finite generating set $X \subset L(\Gamma)$.

Our result settles the question of the value of $\delta_0(X)$ for an 
arbitrary set of self-adjoint generators of a property T factor in full generality: 
\begin{result} Suppose that $N$ is a property T diffuse von Neumann 
algebra 
with a finite set of selfadjoint generators $X$, and let $R^\omega$ be an 
ultrapower of the hyperfinite II$_1$ factor.  Then $\delta_0(X) \leq 1$. 
Moreover, if $N$ has an embedding into $R^\omega$, then $\delta_0(X)=1$,
and if $N$ has no embedding into $R^\omega$, then $\delta_0(X)=-\infty$.
\end{result}

Note that this result shows that the value of 
the free entropy dimension $\delta_0$ is independent of the choice of 
generators of $N$.  In particular, one gets as a corollary that if $\Gamma$ is any 
infinite discrete group with property T, and $X$ is any set of 
self-adjoint generators of the group von Neumann algebra $L(\Gamma)$ (we 
do not make the assumption that $X\subset \mathbb C\Gamma$ here), then 
$\delta_0(X)=1$ or $-\infty$, depending on whether $\Gamma$ embeds into 
the unitary group of $R^\omega$.

The proof of the main theorem relies on a deformation/rigidity argument in 
the style of Popa, which is used to prove that the set of unitary 
conjugacy classes of embeddings of a property T von Neumann algebra $N$ 
into the ultrapower of the hyperfinite II$_1$ factor is discrete.  This 
fact can then be employed to show that if $X\subset N$ is a set of 
self-adjoint generators, then any $k\times k$ matricial microstate for $X$ 
essentially lies in the unitary orbit of a certain discrete set $S$, all 
of whose elements are at least a certain fixed distance apart.  One then 
turns this into an estimate for the packing dimension of the 
microstate space for $X$.  We prove, effectively, that the packing 
dimension of the microstate set is essentially the same as that of a 
small number of disjoint copies of the $k$-dimensional unitary group.

\section{Property T, Embeddings, and Unitary Orbits}

Throughout this section and the next we fix a property T finite von 
Neumann algebra $N$ and a finite $p$-tuple of selfadjoint generators 
$X\subset N$.  $\|\cdot \|_2$ denotes the $L^2$-norm induced by a 
specified trace on a von Neumann algebra.  $M^{sa}_k(\mathbb C)$ denotes 
the set of selfadjoint $k\times k$ matrices and $M_k(\mathbb C)$ denotes 
the set of $k \times k$ matrices.  $tr_k$ is the trace on $M_k(\mathbb 
C)$. If $\xi = \{y_1,\ldots, y_p\}$ and $\eta = \{z_1,\ldots, z_p\}$ are 
$p$-tuples in a von Neumann algebra and $u, w$ are element in a tracial 
von Neumann algebra, then $\xi-\eta = \{y_1 -z_1,\ldots, y_p -z_p\}$, 
$u\xi w= \{uy_1w, \ldots, uy_pw\}$, and $\|\xi \|_2 = (\sum_{i=1}^p 
\|y_j\|_2^2)^{\frac{1}{2}}$. $R >0$ will be a fixed constant greater 
than any of the operator norms of the elements in $X$.   
$\Gamma_R(X;m,k,\gamma)$ will denote the
standard microstate spaces introduced in \cite{dvv:entropy2}. 

The following theorem, stated for the reader's convenience, 
is by now among the standard results in the theory of rigid factors.
Such deformation-conjugacy arguments have played a 
fundamental role in the recent startling results of 
Popa and others
(\cite{ioana-peterson-popa:noOut}, \cite{ozawa},\cite{popa:betti}, 
\cite{popa-sinclair-smith:perturbations}).

\begin{theorem} \label{thm:conjugacy} Let $X$ and $N$ be as above.  Then 
for any $t>0$ there exists a corresponding $r_t >0$ so that if $(M,\tau)$ 
is a tracial von Neumann algebra and $\pi, \sigma: N \rightarrow M$ are 
normal faithful trace-preserving $*$-homomorphisms such that for all $x 
\in X$, $\|\pi(x) -\sigma(x) \|_2 <r_t$, then there exist  
projections $e \in \pi(N)^{\prime} \cap M$, $f \in 
\sigma(N)^{\prime}\cap M$, a partial isometry $v \in M$ such that 
$v^*v = e$, $vv^* =f$, $\tau(e) > 1-t$, and for all $x \in N$,
$v e \pi(x) e v^* = f \sigma(x) f$. \end{theorem}

\begin{proof} Recall (see [15]) that there exist $K, \epsilon_0 >0$, and a 
finite set $F \subset N$ such that if $0 < \delta \leq \epsilon_0$ and $H$ 
is a correspondence of $N$ with a vector $\xi \in H$ satisfying $\|z \xi - 
\xi z\|_2 < \delta$, $z \in F$, then there exists a vector $\eta \in H$ 
which is central for $M$ and $\| \eta - \xi \|_2 < K \delta$. Choose $r_t$ 
so small so that if $\rho_1, \rho_2 : N \rightarrow M$ are any two 
faithful, normal trace preserving $*$-homomorphisms such that for all $x 
\in X$, $\|\rho_1(x) - \rho_2(x)\|_2 < r_t$, then for all $z 
\in F$, $\|\rho_1(z) - 
\rho_2(z)\|_2 < \min\{t, \epsilon_0\} \cdot (4K)^{-1}$.  This can be done 
because $X$ generates $N$.

Suppose $\pi, \sigma: N \rightarrow M$ are two normal, faithful 
trace-preserving $*$-homomorphisms such that for all $x \in X$, $\| \pi(x) 
- \sigma(x) \|_2 < r_t$.  Consider $L^2(M)$ as an $N-N$ bimodule where for 
  any $\xi \in L^2(M)$, $x,y \in N$, $x \xi y = \pi(x)J\sigma(y)^*J \xi$.  
  Denote by $1_M$ the vector associated to the unit of $M$. The 
  hypothesis on $\pi$ and $\sigma$ guarantee that for all $x \in F$, $\|x 
  1_M - 1_M x\|_2 = \|\pi(x) - \sigma(x)\|_2 < \min\{t, \epsilon_0\} 
  \cdot (4K)^{-1}$ which in turn implies the existence of a central 
vector 
  $\eta_0 \in L^2(M)$ for $N$ such that $\|\eta_0 - 1_M \|_2 < t/4$.  
Regard $\eta_0$ as an unbounded operator on $L^2(M)$ by its left action.  
If $\eta_0 = u|\eta_0|$ is the polar decomposition of $\eta_0$, then  $u \in M$ and $\| \eta_0 - 1_M \|_2 < t/4 \Rightarrow \|u - 1_M\|_2 < t/2 \Rightarrow \|u^* u - 1_M \|_2 < t$.  On 
the other hand, since for any $x \in N$, $x \eta_0 = \eta_0 x$, one 
concludes in the usual way that $x u = ux$.  Consequently, $uu^* \in 
\pi(N)^{\prime}$, $u^*u \in \sigma(N)^{\prime}$.  Set $e = uu^* \in 
\pi(N)^{\prime} \cap M$ and $f =u^* u \in \sigma(N)^{\prime} \cap M$.  It 
follows that for all $x \in N$, $u^* e \pi(x) e u = f \sigma(x) f$.  
Finally, $\tau(e) = \tau(f) > 1-t$.\end{proof}

For each $t>0$, we now choose a critical $r=r_t >0$ dependent on $t$ 
as in Theorem~\ref{thm:conjugacy}.

We now need some notation.
\begin{notation} (a) If $\eta \in (M^{sa}_k(\mathbb C))^p$ and $r>0$, then 
$$\Theta_r(\eta) = \{ \xi \in (M^{sa}_k(\mathbb C))^p: \text{for some } 
u \in U_k, \|\xi - u^* \eta u\|_2 < r\}.$$
(b) If $\eta \in (M^{sa}_k(\mathbb C))^p$ and $\kappa, s 
>0$, then $\mathcal G_{\kappa, s}(\eta)$ consist of all all $p$-tuples 
$\xi$ such that there exists projections $e, f \in M^{sa}_k(\mathbb C)$ 
and a $w \in M_k(\mathbb C)$ such that $w^*w =e$, $ww^* =f$, $\tr_k(e) = 
\tr_k(f) > s$ and $\|we\xi ew^* -f\eta f\|_2 < \kappa$. \end{notation}

\begin{lemma}\label{lemma:conj-microstates} For any $\kappa,t  >0$ there exist
an $m \in \mathbb N$ such that if $\xi, \eta \in
\Gamma_R(X;m,k, m^{-1})$ and $\xi \in \Theta_{r_t}(\eta)$, then $\xi \in
\mathcal G_{\kappa, 1-t}(\eta)$. \end{lemma}

\begin{proof} We proceed by contradiction.  Assume that there exists some 
$\kappa_0, t_0 >0$ such that for each $m \in \mathbb N$ there are $k_m \in 
\mathbb N$ and $\xi_m, \eta_m \in \Gamma_R(X;m,k_m, m^{-1})$ such that 
$$\xi 
\in \Theta_{r}(\eta)\textrm{ and }\xi \notin \mathcal G_{\kappa_0, 
1-t_0}(\eta).$$ Fix a free ultrafilter $\omega$, and consider the 
ultraproduct $$R^\omega = \prod^\omega M_{k_m}(\mathbb C) = \frac{ 
\prod_{m=1}^\infty M_{k_m}(\mathbb C) } {\{\langle 
x_m\rangle_{m=1}^\infty: \lim_\omega \\tr_{m_k}(x_m^* x_m) = 0\}}.$$ 
Denote by $Q : \prod M_{k_m} \rightarrow R^{\omega}$ the quotient map.  
Set $\xi = \langle \xi_m \rangle_{m=1}^{\infty}$ and $\eta = \langle 
\eta_m \rangle_{m=1}^{\infty}$.

 For each $m$ we can find a $k_m \times k_m$ unitary $u_m$ such that 
$\|u^*_m \xi_m u_m - \eta\|_2 < r$.  Set $u = \langle u_m 
\rangle_{m=1}^{\infty}$.  It follows that there exist two normal faithful 
trace-preserving $*$-homomorphisms $\pi, \sigma: N \rightarrow R^{\omega}$ 
such that $\pi(X) = Q(U)^* Q(\xi) Q(U)$ and $\sigma(X) = Q(\eta)$.  
Clearly $\|\pi(X) - \sigma(X)\|_2 <r$.  By Theorem~\ref{thm:conjugacy} 
there exist projections $e \in \pi(N)^{\prime} \cap 
R^{\omega}$, $f \in \sigma(N)^{\prime} \cap R^{\omega}$ and a partial 
isometry $v \in R^{\omega}$ with initial domain $e$ and final range 
$f$ such that for all $x \in N$, $ve \pi(x) e v^* 
= f \sigma(x)  f$ and $\tau(e) = \tau(f) > 1- t_0$.
$v$ is a partial isometry 
and $\tau(v^*v) = \tau(e) > 1- t_0$.  There exist sequences of 
projections $\langle e_m \rangle_{m=1}$ and $\langle f_m 
\rangle_{m=1}^{\infty}$ such that for each $m$, $e_m, f_m \in 
M_{k_m}(\mathbb C)$ and $Q(\langle e_m \rangle_{m=1}^{\infty}) = e$, $Q 
(\langle f_m \rangle_{m=1}^{\infty}) =f$.  Similarly there exists a 
sequence of partial isometries $\langle v_m \rangle_{m=1}^{\infty}$ such 
that for each $m$, $v_m \in M_{k_m}(\mathbb C)$ and $Q(\langle v_m 
\rangle_{m=1}^{\infty}) =v$.  We can also arrange it so that for each $m$, 
$v_m v_m^* = f_m$ and $v_m^* v_m = e_m$.  Now, the equation $ve\pi(x)ev^* 
= f \sigma(x) f$, $x\in M$ implies in particular, that for some $\lambda_0 
\in \omega$ $\|v_{m_{\lambda_0}} e_{m_{\lambda_0}} \xi_{m_{\lambda_0}} 
e_{m_{\lambda_0}} v^*_{m_{\lambda_0}} - f_{m_{\lambda_0}} 
\eta_{m_{\lambda}} f_{m_{\lambda_0}} \|_2 < \kappa_0$ and that 
the normalized trace of both $f_{m_{\lambda_0}}$ and $e_{m_{\lambda_0}}$ is strictly 
greater than $1 - t_0$.  But this means that $\xi_{m_{\lambda_0}} \in 
\mathcal G_{\kappa_0, 1-t_0}(\eta)$ which contradicts our initial 
assumption.
\end{proof}

\begin{remark} Observe that in Lemma~\ref{lemma:conj-microstates} the quantity $r_t$ is independent of $\kappa$. 
\end{remark}

\section{The Main Estimate}

In this section we maintain the notation for $\mathbb K_{\epsilon}$ 
introduced in \cite{jung:packing} taken now with respect to the microstate 
spaces with the operator norm cutoffs.  Set $K = \|X \|_2$.  We first 
state a technical lemma on the covering numbers for the spaces $\mathcal 
G_{\kappa, s}(\eta)$.

\begin{lemma} If $\eta \in (M^{sa}_k(\mathbb C))^p$ and $\epsilon, \kappa, 
s >0$ with $\epsilon > \kappa$, then there exists an $5K \epsilon$-net for 
$\mathcal G_{\kappa, s}(\eta)$ with cardinality no greater than
\[  \left ( \frac{2\pi}{\epsilon} \right)^{2k^2- s^2k^2} \cdot
\left(\frac{K+1}{\epsilon} \right)^{4(1-s)^2 k^2}. \]
\end{lemma}

\begin{proof} Find the smallest $m \in \mathbb N$ such that $sk \leq m 
\leq k$.  Denote by $V$ the set of partial isometries in $M_k(\mathbb C)$ 
whose range has dimension $m$.  Denote by $P_m$ the set of projections of trace $mk^{-1}$.  It follows from \cite{szarek:metricentropy} that there exists an $\epsilon$-net for $P_m$ (with respect to the operator norm) with cardinality no greater that $(\frac{2\pi}{\epsilon})^{k^2 - m^2 - (k-m)^2}$.  There exists again by \cite{szarek:metricentropy} an $\epsilon$-net for the unitary group of $M_m(\mathbb C)$ (with respect to the operator norm) with cardinality no greater than $(\frac{2 \pi}{\epsilon})^{m^2}$.  These two facts imply that there exists an $\epsilon$-net $\langle v_{jk} \rangle_{j \in J_k}$ for $V$ with respect to the operator norm such that

\[ \#J_k < \left (\frac{2 \pi}{\epsilon} \right)^{2km - m^2} .
\]
Now fix $j \in J_k$.  Denote by $G(\eta, j)$, the set of all 
$\xi \in (M^{sa}_k(\mathbb C))^p$ such that $\|\xi\|_2 \leq K$ and  
$\|v_{jk}(e_{jk} \xi e_{jk} )v_{jk}^* - f_{jk} 
 \eta f_{jk}\|_2 < 5 K \epsilon $ where $e_{jk} =v_{jk}^*v_{jk}$ 
and $f_{jk}= v_{jk} 
v_{jk}^*$.  There exists a $2\epsilon$-cover $\langle \xi_{ijk} \rangle_{i 
\in \theta(j)}$ for $G(\eta, j)$ such that $\# \theta(j) < 
\left(\frac{K+1}{\epsilon}\right)^{4(1-s)^2k^2}$. 

Consider the set $\langle \xi_{ijk} \rangle_{i \in \theta(j), j \in J_k}$.  
It is clear that this set has cardinality no greater than
\[ \left ( \frac{2\pi}{\epsilon} \right)^{2km - m^2} \cdot 
\left(\frac{K+1}{\epsilon} \right)^{4(1-s)^2 k^2}.  
\]
It remains to show that this set is a $5K\epsilon$-cover for 
$\mathcal G_{\kappa, s}(\eta)$.  Towards this end suppose $\xi \in 
\mathcal G_{\kappa, s}(\eta)$.  Then there exists a partial isometry $v 
\in M_k(\mathbb C)$ such that $v^* v =e$, $vv^* =f$, $\|ve\xi ev^* - f\eta 
f\|_2 <\kappa$, and $\tr_k(e) = \tr_k(f) > s$.  By cutting the domain and 
range of the projection, we can assume that $e$ and $f$ are projections 
onto subspaces of dimension exactly $m$ and we can assume that the 
inequality with tolerance $\kappa$ is preserved.  Obviously then $v \in 
V$, whence there exists a $j_0 \in J_k$ such that $\|v_{j_0k} - 
v\| 
< \epsilon$.  This condition immediately implies that $\|v_{j_0k} e_{j_0k}
- v e\|, \|f_{j_0k} -f\| < 2\epsilon$ and thus
\[ \|v_{j_0k} e_{j_0k} \xi e_{j_0k} v_{j_0k}^* - f_{j_0k} \eta f_{j_0k}\|_2 
\leq 4 \epsilon K + \|ve\xi ev^* - f\eta f\|_2 < 5 K \epsilon.
\]
By definition, $\xi \in G(\eta, j_0)$.  Thus, there exists some 
$i_0$ such that $i_0 \in \theta(j_0)$ and $\|\xi_{i_0j_0k} - \xi\|_2 < 
5K\epsilon$.
\end{proof}

We can now prove the main result of the paper:

\begin{theorem} Let $N$ be a property T diffuse von Neumann algebra with a finite set 
of selfadjoint generators $X$, and let $R^\omega$ be an ultrapower of the 
hyperfinite II$_1$ factor.   \\
(a) If $N$ has an embedding into $R^\omega$, then $\delta_0(X)=1$.
(b) If $N$ has no embedding into $R^\omega$, then $\delta_0(X)=-\infty$.
\end{theorem}

\begin{proof} Fix $1 > a >0$.  For any $\epsilon >0$, setting $\kappa = 
\epsilon$ and $t = 1-a$ in Lemma\ref{lemma:conj-microstates} shows that 
there exists an $m \in \mathbb N$, $m > p^2$, such that if $\xi, \eta \in 
\Gamma_R(X;m,k,m^{-1})$ and $\xi \in \Theta_{r_a}(\eta)$, then $\xi \in 
\mathcal G_{\epsilon, 1-a}(\eta)$.  Consider the ball $B_k$ of 
$(M^{sa}_k(\mathbb C))^p$ of $\| \cdot \|_2$-radius $K+1$.  For each $k$ 
find an $r_a$-net $\langle \eta_{jk} \rangle_{j \in J_k}$ of 
$\Gamma_R(X;m,k,m^{-1})$ with minimal cardinality such that each element 
of the net lies in $\Gamma(X;m,k,m^{-1})$.  The standard volume comparison 
test of this set with $B_k$ (remember that $\Gamma_R(X;m,k,m^{-1}) \subset 
(M^{sa}_k(\mathbb C))^p_K$) implies that \[ \#J_k \leq \left 
(\frac{K+2}{r_a} \right )^{pk^2}. \] For each such $j \in J_k$ find a 
$5K\epsilon$-net $\langle \xi_{ij} \rangle_{i \in \theta(j)}$ for 
$\mathcal G_{\epsilon, 1-a}(\eta_{jk})$ where $\theta(j)$ is an indexing 
set satisfying \[ \# \theta(j) \leq \left ( \frac{2\pi}{\epsilon} 
\right)^{2k^2- (1-a)^2k^2} \cdot \left (\frac{K+2}{\epsilon} 
\right)^{4a^2k^2}. \] Consider now the set $\langle \xi_{ij} \rangle_{i 
\in \theta(j), j \in J_k}$.  It is clear that this set has cardinality no 
greater than \[ \left (\frac{K+2}{r_a} \right )^{pk^2} \left ( 
\frac{2\pi}{\epsilon} \right)^{(1+2a -a^2)k^2} \cdot \left 
(\frac{K+2}{\epsilon} \right)^{4a^2k^2}. \] Moreover, if $\xi \in 
\Gamma_R(X;m,k,m^{-1})$, then there exists some $j_0 \in J_k$ such that 
$\|\xi - \eta_{j_0k}\|_2 < r_a$.  Clearly then, $\xi \in 
\Theta_{r_a}(\eta)$ which implies that $\xi \in \mathcal G_{\epsilon, 
1-a}(\eta_{j_0k})$.  Consequently there exists some $i_0 \in \theta(j_0)$ 
such that $\|\xi - \xi_{i_0 j_0} \|_2 < 5K \epsilon.$ Therefore, $\langle 
\xi_{ij} \rangle_{i \in \theta(j), j \in J_k}$ is a $5K\epsilon$-net for 
$\Gamma_R(X;m,k,m^{-1})$.

The preceding paragraph implies that for $\epsilon >0$,
\begin{eqnarray*} \mathbb K_{5K\epsilon}(X)  & \leq & \limsup_{k 
\rightarrow 
\infty} k^{-2} \cdot \log \left [  \left (\frac{K+2}{r_a} \right )^{pk^2} 
\left ( \frac{2\pi}{\epsilon}
\right)^{(1+2a-a^2)k^2} \cdot
\left (\frac{K+2}{\epsilon} \right)^{4a^2k^2}\right] \\ & = & p |\log 
r_a| + (1+2a- a^2) |\log \epsilon| + \log \left [ (2\pi)^2 (K+2)^{p+4} 
\right].
\end{eqnarray*}
Keeping in mind that $a$ and $\epsilon$ are independent it now 
follows from \cite{jung:packing} 
\begin{eqnarray*} \delta_0(X) & = & \limsup_{\epsilon \rightarrow 0} 
\frac{\mathbb K_{\epsilon}(X)}{|\log \epsilon|} \\ & = & \limsup_{\epsilon 
\rightarrow 0} \frac{\mathbb K_{5K \epsilon}(X) }{|\log \epsilon|} \\ & 
\leq & \limsup_{\epsilon \rightarrow 0} p \cdot \frac{|\log r_a|}{|\log 
\epsilon|} + 1 + 2a - a^2 + 
\frac{\log\left((2\pi)^2(K+2)^{p+4})\right)}{|\log 
\epsilon|} \\ & = & 1 + 2a - a^2.
\end{eqnarray*}
As $1 > a >0$ was arbitrary, $\delta_0(X) \leq 
1$.  The rest of the assertions follow from \cite{jung-freexentropy}.
\end{proof}

\begin{remark}  For $\epsilon >0$ consider the set $X+ \epsilon S = \{x_1
+ \epsilon s_1, \ldots, x_n + \epsilon s_n\}$ where $\{s_1,\ldots, s_n\}$
  is a semicircular family free with respect to $X$. \cite{brown:perturbations} shows that for
  sufficiently small $\epsilon >0$ the von Neumann algebras $M^{\epsilon}$
  generated by $X + \epsilon S$ are not isomorphic to the free group
  factors and yet, if $X^{\prime \prime}$ embeds into the ultraproduct of
  the hyperfinite $\mathrm{II}_1$-factor, then $\chi(X + \epsilon S) >
  -\infty$.  Theorem 2.2 implies that if $X^{\prime \prime}$ embeds into 
the
  ultraproduct of the hyperfinite $\mathrm{II}_1$-factor, then
  $M^{\epsilon}$ cannot have property T.  Also observe that the usual
rigidity/deformation argument shows that for sufficiently small $\epsilon
>0$, there exists a $\mathrm{II}_1$ property T subfactor
$N^{\epsilon}$ of $M^{\epsilon}$.

\end{remark}

\begin{remark} Unfortunately, we were not able to settle the question of 
whether $N$ must be strongly $1$-bounded in the sense of \cite{jung:onebounded}. 
\end{remark}

\noindent{\it Acknowledgments.} The authors would like to thank Adrian 
Ioana, Jesse Peterson, and Sorin Popa for useful conversations.

\bibliographystyle{amsplain}

\providecommand{\bysame}{\leavevmode\hbox to3em{\hrulefill}\thinspace}

\end{document}